\documentclass[12pt]{article}
\usepackage{amsmath, amsthm, amssymb, amsfonts, enumerate}
\begin{document}
\newcommand{\bea}{\begin{eqnarray}}
\newcommand{\ena}{\end{eqnarray}}
\newcommand{\beas}{\begin{eqnarray*}}
\newcommand{\enas}{\end{eqnarray*}}
\newcommand{\beq}{\begin{equation}}
\newcommand{\enq}{\end{equation}}
\newcommand{\ignore}[1]{}
\newcommand{\btheta}{\mbox{\boldmath {$\theta$}}}
\newcommand{\transpose}{{\mbox{\scriptsize\sf T}}}
\newcommand{\qmq}[1]{\quad \mbox{#1} \quad}
\newcommand{\bbox}{\hfill $\Box$}
\newcommand{\Q}{Q}
\newcommand{\nn}{\nonumber}
\newtheorem{theorem}{Theorem}[section]
\newtheorem{corollary}{Corollary}[section]
\newtheorem{conjecture}{Conjecture}[section]
\newtheorem{proposition}{Proposition}[section]
\newtheorem{lemma}{Lemma}[section]
\newtheorem{definition}{Definition}[section]
\newtheorem{example}{Example}[section]
\newtheorem{remark}{Remark}[section]
\newtheorem{case}{Case}[section]
\newtheorem{condition}{Condition}[section]

\title{{\bf\Large Limiting Behavior of High Order Correlations for Simple Random Sampling}}
\author{Christopher Wayne Walker \\ \\ University of Southern California and \\Northrop Grumman Aerospace Systems}
\footnotetext{cwwalker@cwwphd.com}
\footnotetext{
\ignore{Asymptotic enumeration}, 62D05 \ignore{Sampling theory,
sample surveys} } 
\maketitle

\date{}

\begin{abstract}
For $N=1,2,\ldots$, let ${\cal S}_N$ be a simple random sample of
size $n=n_N$ from a population ${\cal A}_N$ of size $N$, where
$0 \le n \le N$. Then with $f_N=n/N$, the sampling fraction, and
${\bf 1}_A$ the inclusion indicator that $A \in {\cal S}_N$, for
any $H \subset {\cal A}_N$ of size $k \ge 0$, the high order correlations
$$
\mbox{Corr}(k) = E \left( \prod_{A \in H} ({\bf 1}_A-f_N) \right)
$$
depend only on $k$, and if the sampling fraction $f_N \rightarrow
f$ as $N \rightarrow \infty$, then 
$$
N^{k/2}\mbox{Corr}(k) \rightarrow [f(f-1)]^{k/2}EZ^k,\ k\ even
$$
and
$$N^{(k+1)/2} \mbox{Corr}(k)\rightarrow[f(f-1)]^{(k-1)/2}(2f-1)\frac{1}{3}(k-1)EZ^{k+1},\ k\ odd$$
where $Z$ is a standard normal random variable. This proves 
a conjecture given in [\ref{AGL}].

\end{abstract}

\section{Introduction}
Simple random sampling is doubtless one of
the most often used tools in statistics [\ref{Hajek}], and it
might appear that nothing new regarding this simple scheme remains
unexplored. With $0 \le n \le N$, by a simple random sample
of size $n$ from a set ${\cal A}_N$ of size $N$ we mean the random subset
of ${\cal S}_N$ of ${\cal A}_N$ with distribution
\bea
\label{Prob=N-choose-n-inv} P({\cal S}_N=r)=\left\{
\begin{array}{cc}
{N \choose n}^{-1} & r \subset {\cal A}_N, |r|=n\\
0                  & \mbox{otherwise,}
\end{array}
\right.
\ena
where $|r|$ denotes the size of the set $r$.
It is easy to see that all individuals in ${\cal A}_N$
have an equal chance of being included in the sample,
and that in particular, for $A \in
{\cal A}_N$, \bea \label{def-fN} E {\bf 1}_A =f_N \quad
\mbox{where} \quad f_N=\frac{n}{N}, \ena where the inclusion
indicator ${\bf 1}_A$ takes the value 1 when $A \in {\cal S}_N$
and the value 0 otherwise. The value $f_N$ is known as the
sampling fraction. Likewise, from (\ref{Prob=N-choose-n-inv}) one can
show directly that the inclusion
indicators ${\bf 1}_A$ and ${\bf 1}_B$, for distinct individuals
$A$ and $B$ in ${\cal A}_N$, are negatively correlated, as the
inclusion of $A$ leaves less room in the remaining sample for $B$.
That is, \beas
E{\bf 1}_A {\bf 1}_B = \frac{n(n-1)}{N(N-1)} < \left( \frac{n}{N}\right)^2 = E{\bf 1}_A
E{\bf 1}_B,
\enas
or, considering the two way correlation \bea
\label{Corr2} \mbox{Corr}({\bf 1}_A,{\bf 1}_B)=E\left( ({\bf
1}_A-f_N)({\bf 1}_B-f_N\right), \ena we have
\bea
\label{Corr2SRS} \mbox{Corr}({\bf 1}_A,{\bf 1}_B)=
\frac{-n(N-n)}{N^2(N-1)}. \ena However, the higher order
correlations of simple random sampling, which arise in some
applications [\ref{AGL}] and exhibit rather interesting behavior,
are virtually unknown.

To consider such correlations of higher order,
generalizing (\ref{Corr2}), for any $H \subset {\cal
A}_N$ of size $|H|=k, 0 \le k \le n$, we define
\bea \label{def-Corr-H}
\mbox{Corr}(H)=E \left( \prod_{A \in H}({\bf 1}_A-f_N) \right).
\ena
We see from (\ref{Prob=N-choose-n-inv}) that the probability that
all individuals in $H$ are included in the sample is \bea
\label{E-prod-1A} E\left( \prod_{A \in H}{\bf 1}_A
\right)=\frac{{N-k \choose n-k}}{{N \choose n}}, \ena which we
note only depends on $k$, and not on which individuals comprise
the set $H$. Hence, arguing either directly using
(\ref{E-prod-1A}), or by noting the indicators $\{{\bf 1}_A, A \in
{\cal A}_N\}$ are exchangeable, we conclude that $\mbox{Corr}(H)$
depends only on the size $k$ of $H$, and hence we denote it by
$\mbox{Corr}(k)$.

In [\ref{AGL}] the high order correlation of rejective sampling
was studied in order to determine the asymptotic properties of a
generalized logistic estimator. In rejective sampling
[\ref{Hajek2}], each individual $A$ in ${\cal A}_N$ is associated
with a non-negative weight $w_A$, and the probability of sampling
a set $r \subset {\cal A}_N$ of size $n$ is given by
$$
P({\cal S}_n=r) = \frac{w_r}{\sum_{s \subset {\cal A}_N,
|s|=n}w_s} \quad \mbox{where} \quad w_s=\prod_{j \in s}w_j.
$$
We note that the high order correlations of rejective sampling may
be defined exactly as in (\ref{def-Corr-H}), and that simple
random sampling is the special case of rejective sampling, taking
all weights equal.

Critical in the asymptotic analysis in [\ref{AGL}] was the fact
that under some stability conditions on the weights, the second
and third order correlations of rejective sampling decay at rates
$N^{-1}$ and $N^{-2}$, respectively, that is, that
$$
\lim_{N \rightarrow \infty} N\mbox{Corr}(2)=O(1), \quad \mbox{and}
\quad \lim_{N \rightarrow \infty} N^2\mbox{Corr}(3) = O(1).
$$

Checking for the special case of simple random sampling when $f_N
\rightarrow f$ as $N \rightarrow \infty$, equality (\ref{Corr2SRS}) implies
$N\mbox{Corr}(2) \rightarrow f(1-f)$, and further direct
calculation for correlations up to order 9 obtained by expanding
the expression in definition (\ref{def-Corr-H}) and using (\ref{E-prod-1A}) yields
\beas
N   \mbox{Corr}(2) &\rightarrow& f(f-1)\\
N^2 \mbox{Corr}(3) &\rightarrow& 2f(f-1)(2f-1)\\
N^2 \mbox{Corr}(4) &\rightarrow& 3f^2(f-1)^2\\
N^3 \mbox{Corr}(5) &\rightarrow& 20f^2(f - 1)^2(2f - 1)\\
N^3 \mbox{Corr}(6) &\rightarrow& 15f^3(f - 1)^3\\
N^4 \mbox{Corr}(7) &\rightarrow& 210f^3(f - 1)^3(2f-1)\\
N^4 \mbox{Corr}(8) &\rightarrow& 105f^4(f - 1)^4\\
N^5 \mbox{Corr}(9) &\rightarrow& 2520f^4(f - 1)^4(2f-1). \enas
Perhaps the most surprising feature of these correlations is that
their rate of decay depends on the parity of the correlation
order, in particular, one can easily conjecture that \bea
\label{parity} N^{(k+k\,mod \,2)/2} \mbox{Corr}(k) = O(1) \quad
k=1,2,\ldots. \ena Theorem 4.1 of [\ref{AGL}] shows that
(\ref{parity}) holds quite generally for rejective sampling, and
therefore for simple random sampling in particular. Application of
this theorem sufficed to complete the asymptotic analysis required
in [\ref{AGL}] for rejective sampling.

Another feature of the simple random sampling
correlations is also easy to conjecture, that their scaled limits
are equal to a constant depending on $k$, times the factor
$f(f-1)$ raised to $(k-k\,mod \,2)/2$, and if $k$ is odd, times
the additional factor $(2f-1)$. Hence, one need only determine the
constants to completely specify their limiting behavior. On the
basis of the above observations and the constants corresponding to
the even and odd values of $k$ up to 9, that is, the sequences
$1,3,15,105$ and $2,20,210,2520$, respectively, a conjecture was
put forth in [\ref{AGL}], which is now validated by the following
theorem which is proven in this analysis.

\begin{theorem}
\label{main}
For $N=1,2,\ldots$ let ${\cal S}_N$ be a sequence of simple random samples
from populations ${\cal A}_N$ of size $N$, whose sampling fractions $f_N$ obey
$$
f_N \rightarrow f \in (0,1) \quad \mbox{as $N \rightarrow
\infty$}.
$$
Then
$$\lim_{N\rightarrow \infty} N^{k/2} \mbox{Corr}(k)=[f(f-1)]^{k/2} EZ^k,\ k\ even$$
and
$$\lim_{N\rightarrow \infty} N^{(k+1)/2} \mbox{Corr}(k)=[f(f-1)]^{(k-1)/2}(2f-1)\frac{1}{3}(k-1)EZ^{k+1},\ k\ odd$$
where $Z$ is a standard normal variate.
\end{theorem}

We recall that the standard normal variable $Z$ is the one with
distribution function
$$
P(Z \le x) = \frac{1}{\sqrt {2 \pi}} \int_{-\infty}^x e^{-u^2/2}
du,
$$
or, equivalently, moment generating function $M(t)=Ee^{tZ}$ given by
$$
M(t)=e^{t^2/2}.
$$
Substituting $t^2/2$ for $u$ in the expansion
$e^u=\sum_{j=0}^\infty u^j/j!$, we find
$$
M(t)=\sum_{j=0}^\infty \frac{(t^2/2)^j}{j!} = \sum_{j=0}^\infty
\frac{t^{2j}}{2^j j!}=\sum_{j=0}^\infty \frac{(2j)!}{2^j
j!}\frac{t^{2j}}{(2j)!}=\sum_{k=0,2,4,\ldots}
\frac{k!}{2^{k/2}(k/2)!}\frac{t^k}{k!},
$$
and hence, as $EZ^k=M^{(k)}(0)$, we conclude that $Z$ has vanishing
moments of odd order, while
$$
EZ^k =\frac{k!}{2^{k/2}(k/2)!} = \prod_{0< j < k, j \mbox{{\tiny
odd}}}j = \frac{k!}{k(k-2)\cdots 2} \quad \mbox{for
$k$ even}.
$$

The appearance of the moments of the standard normal variable in
counting contexts is well known. In particular, in the On-Line
Encyclopedia of Integer Sequences [\ref{oeis}] there are roughly
twenty combinatorial structures listed which are counted by the
sequence $EZ^{k}, k=0,2,4,\ldots$, such as the number of perfect
matchings in the complete graph $K(2n)$, and the number of
fixed-point-free involutions in the symmetric group $S_{2n}$.
However, in contrast, the sequence $(1/3)(k-1)EZ^{k+1},\ k=1,3,5,\ldots$ has only
two listings, and of a more abstract nature, in particular, the
Ramanujan polynomials $-\psi_{n+2}(n+2,x)$ evaluated at 1, and,
with offset 2, the second Eulerian transform of
$0,1,2,3,4,\ldots$. In addition to providing a concrete situation
in which this lesser known second sequence appears, the present
work also demonstrates that there exists a connection between it
and the one much better known. At the same time, we remark that
our proof of Theorem \ref{main} is purely computational, and that
it would be enlightening to also complete the argument in a manner
which makes the appearance of the normal moments more natural.

\section{Main Result}
The proof of Theorem \ref{main} requires a few identities for the
Stirling and Bernoulli numbers, which can be found in
[\ref{Concrete}]. Letting $(x)_j$ denote the falling factorial, or Pochhammer symbol,
$$
(x)_j=x(x-1)\cdots(x-j+1),
$$
expanding $(x)_j$ as a polynomial in $x$ we have
\beas
(x)_j = \sum_{v=0}^j \left[j \atop v \right]
(-1)^{j-v} x^v,
\enas where $\left[j \atop v \right]$ is the
(unsigned) Stirling numbers of the first kind; for instance,
\beas
\left[j \atop j \right]=1, \quad \left[j
\atop j-1 \right]= \frac{1}{2}j(j-1), \quad \left[j
\atop 1 \right]=(j-1)! \quad \mbox{and} \quad \left[j \atop 0 \right]=0.
\enas
Regarding Stirling numbers of the second kind, denoted $\left\{m \atop k \right\}$, we make use
of the identity
\bea \label{Stirling 2} \sum_{j=0}^k {k\choose
j}(-1)^j j^m=(-1)^k k! \left\{m \atop k \right\}.
\ena
In particular we note that
\bea \label{second_Stirlings} \left\{m \atop
k \right\}=0\ \mbox{for }m<k,\quad \left\{k \atop k \right\}=1,
\quad \mbox{and} \quad \left\{k+1 \atop k \right\}=\left(k+1 \atop
2 \right).\ena Furthermore,  \bea \label{Bernoulli-1}
\sum_{p=0}^{k-1} p^m=\frac{1}{m+1}\sum_{p=0}^m {m+1 \choose p} B_p
k^{m+1-p} \ena where $B_p$ denotes the $p^{th}$ Bernoulli number
defined by the explicit recurrence \bea
\label{Bernoulli-2}\sum_{j=0}^m {m+1 \choose j}B_j=\delta(m), \ena
where $\delta(m)$ is the Kronecker delta function, taking the
value 1 for $m=0$ and zero otherwise.

For any integer
$m$ let $r_{\boldsymbol{\alpha},m}(x)$ denote a polynomial in $x$ of degree no more than $\max(m,0)$ whose coefficients are a function of $\boldsymbol{\alpha}$, possibly a vector; the polynomial $r_{\boldsymbol{\alpha},m}(x)$ is not necessarily the same at each occurrence. Note with this notation that $r_{\boldsymbol{\alpha},m}(x)+r_{\boldsymbol{\alpha},m}(x)=r_{\boldsymbol{\alpha},m}(x)$ and if $f(\boldsymbol{\alpha})$ is a function of $\boldsymbol{\alpha}$ not involving $x$ then $f(\boldsymbol{\alpha})r_{\boldsymbol{\alpha},m}(x)=r_{\boldsymbol{\alpha},m}(x)$. From (\ref{Bernoulli-2}) one easily finds that
$B_0=1$ and $B_1=-1/2$, and hence (\ref{Bernoulli-1}) yields
\bea
\label{expansioncoefficients}
\sum_{p=0}^{k-1}
p^m=\frac{1}{m+1}k^{m+1}-\frac{1}{2}k^{m}+r_{m,m-1}(k).
\ena
Other examples that illustrate this notation are
\beq
\nonumber k^4-2km+m^2+j^2/m-2j=r_{(k,m),2}(j)
\enq
and
\beq
\nn k^4-2km+m=r_{(k,m),0}(j).
\enq

Similarly, let $r_{\boldsymbol{\alpha},m,n,p}(x,y)$ denote a polynomial in $x$ and $y$ whose coefficients are a function of $\boldsymbol{\alpha}$ where the highest power of $x$ that occurs is at most $\max(m,0)$, the highest power of $y$ that occurs is at most $\max(n,0)$, and the
highest sum $a+b$ in terms of the form $x^ay^b$ is at most $\max(p,0)$; again the polynomial $r_{\boldsymbol{\alpha},m,n,p}(x,y)$ is not necessarily the same at each occurrence. Examples that illustrate this notation are
\beq
\nonumber k^4+k^3+k^2 j^2+j^3+j^4=r_{1,4,4,4}(j,k)
\enq
and
\bea
\nonumber k^4+mk^4 j^2+k^2j^2+j^6&=&r_{m,6,4,6}(j,k)\\
\nn &=&r_{(k,m),6}(j).
\ena

An identity that can be obtained directly from Equation 5.114 in [\ref{Concrete}] is
\beq
\sum_{n=0}^{m-1} (-1)^n {m \choose n}=(-1)^{m-1}
\enq
for $m\geq 1$. We can manipulate this last expression 
to deduce
\beq \label{unity identity}
\sum_{n=0}^{m-1} (-1)^n {m \choose n+1}=u(m-1)
\enq
where $u(n)$ is the unit step function defined to be 1 when $n\geq 0$ and is zero otherwise.

Using the well known identity
\bea
{{m}\choose n+1}={{m-1}\choose n+1}+{{m-1}\choose n}
\ena
and (\ref{unity identity}) we easily find 
\bea \label{sum_identity_1}
\sum_{n=0}^{m-1} (-1)^n {m-1\choose n}=\delta(m-1).
\ena


The following can be found in [\ref{Garrappa}]:
\bea \label{bin_sum_identity}
\sum_{n=0}^{m-1} (-1)^n {{m-1}\choose n}\frac{1}{n+\beta}=\frac{\Gamma(m) \Gamma \left(\beta\right)}{\Gamma \left(m+\beta\right)}u(m-1)
\ena
where, $\Gamma(m)$ denotes the gamma function defined by
\beq
\Gamma(z)=\int_0^\infty t^{z-1}e^{-t}dt,\ \ (\Re{z}>0).
\nn 
\enq

We note from [\ref{AS}] that for positive integers $m$,
\beq
\label{Gamma1} \Gamma(m)=(m-1)!,\ \ \Gamma\left(\frac{1}{2}\right)=\sqrt{\pi},\ \ \Gamma\left(\frac{3}{2}\right)=\frac{\sqrt{\pi}}{2}, 
\enq
and
\beq
\nn \Gamma \left(m+\frac{1}{2}\right)=\frac{(2m-1)!!}{2^m}\sqrt{\pi}
\enq
where !! denotes the double factorial and when $m=0$ we define $(-1)!!=1$. We can write the last result as
\bea
\label{Gamma2} \Gamma \left(m+\frac{1}{2}\right)&=&\frac{(2m-1)!}{2^{2m-1}(m-1)!}\sqrt{\pi}.
\ena

Using (\ref{sum_identity_1}) and (\ref{bin_sum_identity}) we get
\bea
\label{sum_identity_3} \sum_{n=0}^{m-1} (-1)^n {{m-1}\choose n}\frac{n}{n+\beta}&=&\delta(m-1)-\frac{\beta \Gamma(m) \Gamma \left(\beta\right)}{\Gamma \left(m+\beta\right)}u(m-1).
\ena

Combining (\ref{bin_sum_identity}) and (\ref{sum_identity_3}) we obtain for real numbers $\alpha, \delta, \gamma, \beta$ and integer $m$ such that $\gamma\neq 0$ and
$\frac{\beta}{\gamma}\neq 0,-1,-2,\ldots,-(m-1)$,
\bea \label{bin_sum_lemma_result}
\lefteqn{\sum_{n=0}^{m-1} (-1)^n {{m-1}\choose n}\frac{\alpha n+\delta}{\gamma n+\beta}}\nn \\
\label{bin_sum_lemma_new} &=&\frac{\alpha}{\gamma}\delta(m-1)+\frac{\Gamma(m)\Gamma \left(\frac{\beta}{\gamma} \right)}{\Gamma \left(m+\frac{\beta}{\gamma} \right)}\frac{1}{\gamma}\left(\delta-\frac{\alpha \beta}{\gamma}\right)u(m-1).
\ena

We now derive another result that will be useful. Using (\ref{expansioncoefficients}) and the binomial theorem we can write for integers $n,p,s\geq0$,
\bea
\lefteqn{\sum_{q=1}^{k-m} q(q+1)^{pn+s}}\nn \\
\nn &=&\sum_{q=1}^{k-m} \sum_{j=0}^{pn+s}{pn+s\choose j}q^{pn+s+1-j}\\
\nn &=&\sum_{q=1}^{k-m} \left(q^{pn+s+1}+(pn+s)q^{pn+s}+r_{(n,p,s),pn+s-1}(q)\right)\\
\label{binomial_Pkm_result} &=&\frac{k^{pn+s+2}}{pn+s+2}\\
\nn &&+\frac{\left((pn+s)(3-2m)+1-2m\right)k^{pn+s+1}}{2(pn+s+1)}\\
\nn &&+r_{(n,m,p,s),pn+s}(k).
\ena
If we let $k=j$ in this last identity and set $m=1$ we get
\bea
\sum_{q=1}^{j-1} q(q+1)^{pn+s}
\label{binomial_Pjn_result} &=&\frac{j^{pn+s+2}}{pn+s+2}+\frac{(pn+s-1)j^{pn+s+1}}{2(pn+s+1)}\\
\nn &&+r_{(n,p,s),pn+s}(j).
\ena

For the following definition we adopt the empty sum convention that
\beas
\sum_{q=a}^b x_q=0 \quad \mbox{when $b<a$.}
\enas
\begin{definition}
\label{P_def}
For nonnegative integers $k,j$ and $m$, let
\beas
\lefteqn{P_{k,0}(j)=1,}\\
&&  P_{k,m}(j)=\sum_{q=1}^{k-m} q P_{k,m-1}(q+1)-\sum_{q=1}^{j-1} q
P_{k,m-1}(q+1) \quad \mbox{for $m \ge 1$,}
\enas
and
\beas
P_{k,0}^0(j)=1 \qmq{and} P_{k,m}^0(j)= \sum_{q=j}^{k-m} q P_{k,m-1}^0(q+1) \quad \mbox{for $m \ge 1$.}
\enas
\end{definition}

\noindent For example, since $\sum_{q=1}^{k-1} q = k(k-1)/2$ for all $k \ge 0$,
for $m=1$ we obtain
\beas
P_{k,1}(j)=\frac{k(k-1)}{2}-\frac{j(j-1)}{2} \quad \mbox{and}
\quad P_{k,1}^0(j)= \left\{
\begin{array}{cc}
P_{k,1}(j) & 0 \le j \le k-1\\
0 & j \ge k.\\
\end{array}
\right.
\enas
The following lemma shows that the identity $P_{k,1}(k)=0$ is a specific instance of a more general fact.

\begin{lemma} \label{P-lemma}
For all nonnegative integers $k,j$ and $m$,
\bea
\label{P2b}
P_{k,m}(j)=0 \quad \mbox{for $0 \le k-m+1 \le j \le k,$}
\ena
and
\bea \label{P2x}
P_{k,m}(j)=P_{k,m}^0(j) \qmq{for} 0 \le j \leq k.
\ena

Furthermore, for all integers $k,m$ and $j$ satisfying $0
\le m,j \le k$,
\bea  \label{P-km-1}
P_{k,m}(j)=(-1)^m\left(\frac{j^{2m}+\frac{1}{3}m(2m-5)j^{2m-1}}{2^m m!}\right)+r_{(k,m),2m-2}(j),
\ena
and for integers $j$ and $v$ satisfying $0\leq v\leq j$,
\bea \label{P_{j,v}}
P_{j,v}(1)=\frac{j^{2v}-\frac{1}{3}v(2v+1)j^{2v-1}}{2^v v!}+r_{v,2v-2}(j).
\ena

\end{lemma}

\noindent {\bf Proof:} First note that (\ref{P2b}) implies (\ref{P2x}), as the equality holds for $0 \le j \le k-m$ by construction. To argue by induction, we have that (\ref{P2b})
holds when $m=0$, as in this case the premise $0 \le k-m+1\leq j \leq k$ is vacuous. Now assume (\ref{P2b}) holds
for $m \ge 0$. From Definition \ref{P_def}, for $k-m \le j \le k$,
\beas
P_{k,m+1}(j) &=&\sum_{q=1}^{k-m-1} q P_{k,m}(q+1) - \sum_{q=1}^{j-1} q P_{k,m}(q+1)\\
&=&-\sum_{q=k-m}^{j-1} q P_{k,m}(q+1)\\
&=&-\sum_{q=k-m+1}^j (q-1) P_{k,m}(q).
\enas
For $k-m+1 \le q \le j \le k$ we have
$P_{k,m}(q)=0$ by (\ref{P2b}). Hence,
\beas
P_{k,m+1}(j)=0 \quad \mbox{for $k-m \le j \le k$,}
\enas
which is (\ref{P2b}) with $m+1$ replacing $m$, completing the proof of the first two claims.

To prove the rest of the lemma we will show that

\bea P_{k,m}(j) \label{P_is_P_tilde_plus_r}
&=&\sum_{n=0}^m \frac{(-1)^n j^{2n}k^{2m-2n}}{2^m (m-n)! n!}\\
\nn &&+\sum_{n=0}^{m-1} \frac{\frac{1}{3}(2n-3)(-1)^{n+1} j^{2n+1} k^{2m-2-2n}}{2^m(m-1-n)!n!} \\
\nn &&+\sum_{n=0}^{m-1} \frac{\frac{1}{3}(2n-2m-1)(-1)^n  j^{2n} k^{2m-1-2n}}{2^m(m-1-n)!n!}\\
\nn &&+r_{m,2m-2,2m-2,2m-2}(j,k)
\ena
from which both (\ref{P-km-1}) and (\ref{P_{j,v}}) will easily follow. 

Now let
\beq
D_{n,m}=2^m (m-n)! n!.
\enq
We will use induction on $m$ in our proof. Clearly (\ref{P_is_P_tilde_plus_r}) is true for $m=0$. Now assume (\ref{P_is_P_tilde_plus_r}) is true for some $m-1 \geq 0$. Then, to evaluate
\bea
\label{Pkm_recursion}
P_{k,m}(j)=\sum_{q=1}^{k-m} q P_{k,m-1}(q+1)-\sum_{q=1}^{j-1} q P_{k,m-1}(q+1)
\ena
let
\bea 
\nn P_{k,m;1}(j)&=&\sum_{n=0}^{m} \frac{(-1)^n j^{2n}k^{2m-2n}}{2^m (m-n)! n!}\\
\nn P_{k,m;2}(j)&=&\sum_{n=0}^{m-1} \frac{\frac{1}{3}(2n-3)(-1)^{n+1} j^{2n+1} k^{2m-2-2n}}{2^m(m-1-n)!n!}  \\
\nn P_{k,m;3}(j)&=&\sum_{n=0}^{m-1} \frac{\frac{1}{3}(2n-2m-1)(-1)^n  j^{2n} k^{2m-1-2n}}{2^m(m-1-n)!n!} \\
\nn P_{k,m;4}(j)&=&r_{m,2m-2,2m-2,2m-2}(j,k)
\ena
so that (\ref{Pkm_recursion}) becomes
\beq \label{Pkm_parts}
P_{k,m}(j)=\sum_{q=1}^{k-m} q \sum_{t=1}^4 P_{k,m-1;t}(q+1)-\sum_{q=1}^{j-1} q \sum_{t=1}^4 P_{k,m-1;t}(q+1).
\enq
We will evaluate the terms in (\ref{Pkm_parts}) separately. We find for $t=1$
\bea
\nn \lefteqn{\sum_{q=1}^{k-m} q P_{k,m-1;1}(q+1)-\sum_{q=1}^{j-1} q P_{k,m-1;1}(q+1)} \\
\nn &=&\sum_{n=0}^{m-1} \frac{(-1)^n k^{2m-2-2n}}{D_{n,m-1}} \sum_{q=1}^{k-m} q(q+1)^{2n} \\
\label{Pk1_part} &&-\sum_{n=0}^{m-1} \frac{(-1)^n k^{2m-2-2n}}{D_{n,m-1}} \sum_{q=1}^{j-1} q(q+1)^{2n}.
\ena

Using (\ref{binomial_Pkm_result}) and (\ref{binomial_Pjn_result}) with $p=2$ and $s=0$, (\ref{Pk1_part}) becomes
\bea
\nn \lefteqn{\sum_{q=1}^{k-m} q P_{k,m-1;1}(q+1)-\sum_{q=1}^{j-1} q P_{k,m-1;1}(q+1)} \\
\nn &=&\frac{k^{2m}}{2^m}\sum_{n=0}^{m-1} \frac{(-1)^n}{(m-1-n)!(n+1)!}+\frac{k^{2m-1}}{2^m}\sum_{n=0}^{m-1}\frac{(-1)^n(2n(3-2m)+1-2m)}{(2n+1)(m-1-n)!n!}\\ \nn
\nn &&-\sum_{n=0}^{m-1} \frac{(-1)^n j^{2n+2} k^{2m-2-2n}}{2^m (m-1-n)!(n+1)!}-\sum_{n=0}^{m-1} \frac{(-1)^n (2n-1)j^{2n+1} k^{2m-2-2n}}{2^m (m-1-n)!n!(2n+1)}\\
\nn &&+r_{m,2m-2,2m-2,2m-2}(j,k)\\
\nn &=&\frac{k^{2m}}{2^m m!}\sum_{n=0}^{m-1} (-1)^n {m \choose {n+1}}\\
\nn &&+\frac{k^{2m-1}}{2^m (m-1)!}\sum_{n=0}^{m-1} (-1)^n{m-1 \choose n}\frac{2(3-2m)n+1-2m}{2n+1}\\ 
\nn &&+\sum_{n=1}^{m} \frac{(-1)^n j^{2n} k^{2m-2n}}{2^m (m-n)!n!}-\sum_{n=0}^{m-1} \frac{(-1)^n j^{2n+1} k^{2m-2-2n}}{2^m (m-1-n)!n!} \frac{(2n-1)}{(2n+1)}\\
\nn &&+r_{m,2m-2,2m-2,2m-2}(j,k).
\ena
Using (\ref{unity identity}) and (\ref{bin_sum_lemma_result}) with $\alpha=2(3-2m),\ \delta=1-2m,\ \gamma=2,$ and $\beta=1$ this becomes
\bea
\nn \lefteqn{\sum_{q=1}^{k-m} q P_{k,m-1;1}(q+1)-\sum_{q=1}^{j-1} q P_{k,m-1;1}(q+1)} \\
\nn &=&\frac{k^{2m}}{2^m m!}u(m-1)-\frac{k^{2m-1}}{2^{m} (m-1)!}\left(\frac{\Gamma(m)\Gamma \left(\frac{1}{2}\right)}{\Gamma\left(m+\frac{1}{2}\right)}u(m-1)-(3-2m)\delta(m-1)\right)\\
\nn &&+\sum_{n=1}^{m} \frac{(-1)^n j^{2n} k^{2m-2n}}{2^m (m-n)!n!}-\sum_{n=0}^{m-1} \frac{(-1)^n j^{2n+1} k^{2m-2-2n}}{2^m (m-1-n)!n!} \frac{(2n-1)}{(2n+1)}\\
\nn &&+r_{m,2m-2,2m-2,2m-2}(j,k).
\ena
Using (\ref{Gamma1}) and (\ref{Gamma2}) we get
\bea
\nn \lefteqn{\sum_{q=1}^{k-m} q P_{k,m-1;1}(q+1)-\sum_{q=1}^{j-1} q P_{k,m-1;1}(q+1)} \\
\label{Pk1_result} &=&\sum_{n=0}^{m} \frac{(-1)^n j^{2n} k^{2m-2n}}{2^m (m-n)!n!}u(m-1)\\
\nn &&-\frac{2^{m-1}(m-1)! k^{2m-1}}{(2m-1)!}u(m-1)+\frac{k^{2m-1}}{2^m (m-1)!}(3-2m)\delta(m-1)\\ 
\nn &&-\sum_{n=0}^{m-1} \frac{(-1)^n j^{2n+1} k^{2m-2-2n}}{2^m (m-1-n)!n!} \frac{(2n-1)}{(2n+1)}+r_{m,2m-2,2m-2,2m-2}(j,k).
\ena

We next find from (\ref{Pkm_parts}) with $t=2$
\bea
\nn \lefteqn{\sum_{q=1}^{k-m} q P_{k,m-1;2}(q+1)-\sum_{q=1}^{j-1} q P_{k,m-1;2}(q+1)} \\
\nn &=&\sum_{n=0}^{m-2} \sum_{q=1}^{k-m} \frac{\frac{1}{3}(2n-3)(-1)^{n+1} q(q+1)^{2n+1} k^{2m-4-2n}}{2^{m-1}(m-2-n)!n!} \\
\nn &&-\sum_{n=0}^{m-2} \sum_{q=1}^{j-1} \frac{\frac{1}{3}(2n-3)(-1)^{n+1} q(q+1)^{2n+1} k^{2m-4-2n}}{2^{m-1}(m-2-n)!n!}.
\ena
Using (\ref{binomial_Pkm_result}) and (\ref{binomial_Pjn_result}) with $p=2$ and $s=1$ we get
\bea
\nn \lefteqn{\sum_{q=1}^{k-m} q P_{k,m-1;2}(q+1)-\sum_{q=1}^{j-1} q P_{k,m-1;2}(q+1)} \\
\nn &=&\sum_{n=0}^{m-2} \frac{\frac{1}{3}(2n-3)(-1)^{n+1} k^{2m-1}}{2^{m-1}(m-2-n)!n!(2n+3)}+r_{m,2m-2}(k) \\
\nn &&-\sum_{n=0}^{m-2} \frac{\frac{1}{3}(2n-3)(-1)^{n+1} j^{2n+3} k^{2m-4-2n}}{2^{m-1}(m-2-n)!n!(2n+3)}\\
\nn &&+r_{m,2m-2,2m-4,2m-2}(j,k)\\
\nn &=&-\frac{\frac{1}{3} k^{2m-1}}{2^{m-1}(m-2)!}\sum_{n=0}^{m-2} (-1)^n {m-2\choose n}\frac{2n-3}{2n+3}\\
\nn &&-\sum_{n=0}^{m-2} \frac{\frac{1}{3}(2n-3)(-1)^{n+1} j^{2n+3} k^{2m-4-2n}}{2^{m-1}(m-2-n)!n!(2n+3)}\\
\nn &&+r_{m,2m-2,2m-2,2m-2}(j,k). \nn
\ena
Using (\ref{bin_sum_lemma_new}) with $m=m-1,\ \alpha=2,\ \delta=-3,\ \gamma=2,\ \beta=3$ and then using (\ref{Gamma1}) and (\ref{Gamma2}) we obtain
%
\bea
\nn \lefteqn{\sum_{q=1}^{k-m} q P_{k,m-1;2}(q+1)-\sum_{q=1}^{j-1} q P_{k,m-1;2}(q+1)} \\
\label{Pk2_result} &=&\frac{(m-1)!2^{m-1}k^{2m-1}}{(2m-1)!}u(m-2)-\frac{\frac{1}{3} k^{2m-1}}{2^{m-1}(m-2)!}\delta(m-2)\\
\nn &&-\sum_{n=0}^{m-2} \frac{\frac{1}{3}(2n-3)(-1)^{n+1} j^{2n+3} k^{2m-4-2n}}{2^{m-1}(m-2-n)!n!(2n+3)}\\
\nn &&+r_{m,2m-2,2m-2,2m-2}(j,k).
\ena

Next we find from (\ref{Pkm_parts}) with $t=3$
\bea
\nn \lefteqn{\sum_{q=1}^{k-m} q P_{k,m-1;3}(q+1)-\sum_{q=1}^{j-1} q P_{k,m-1;3}(q+1)} \\
\nn &=&\sum_{q=1}^{k-m} \sum_{n=0}^{m-2} \frac{\frac{1}{3}(2n-2m+1)(-1)^n q (q+1)^{2n} k^{2m-3-2n}}{2^{m-1}(m-2-n)!n!} \\
\nn &&-\sum_{q=1}^{j-1} \sum_{n=0}^{m-2} \frac{\frac{1}{3}(2n-2m+1)(-1)^n q (q+1)^{2n} k^{2m-3-2n}}{2^{m-1}(m-2-n)!n!} .
\ena
Using (\ref{binomial_Pkm_result}) and (\ref{binomial_Pjn_result}) with $p=2$ and $s=0$ this becomes
\bea
\nn \lefteqn{\sum_{q=1}^{k-m} q P_{k,m-1;3}(q+1)-\sum_{q=1}^{j-1} q P_{k,m-1;3}(q+1)} \\
\nn &=&\sum_{n=0}^{m-2}\frac{\frac{1}{3}(2n-2m+1)(-1)^n  k^{2m-1}}{2^{m-1}(m-2-n)!n!(2n+2)}+r_{m,2m-2}(k)\\
\nn &&-\sum_{n=0}^{m-2}\frac{\frac{1}{3}(2n-2m+1)(-1)^n j^{2n+2}k^{2m-3-2n}}{2^{m-1}(m-2-n)!n!(2n+2)}+r_{m,2m-3,2m-3,2m-2}(j,k) \\ \nn
\nn &=&\frac{\frac{1}{3} k^{2m-1}}{2^m (m-2)!}\sum_{n=0}^{m-2} (-1)^n {m-2\choose n} \frac{2n+1-2m}{n+1}\\
\nn &&+\sum_{n=1}^{m-1}\frac{\frac{1}{3}(2n-2m-1)(-1)^n j^{2n}k^{2m-1-2n}}{2^{m}(m-1-n)!n!}+r_{m,2m-3,2m-2,2m-2}(j,k).\\ \nn
\ena
Using (\ref{bin_sum_lemma_new}) with $m$ replaced by $m-1$,\ $\alpha=2,\ \delta=1-2m,\ \gamma=1$ and $\beta=1$ this becomes
\bea
\nn \lefteqn{\sum_{q=1}^{k-m} q P_{k,m-1;3}(q+1)-\sum_{q=1}^{j-1} q P_{k,m-1;3}(q+1)} \\
\nn &=&\frac{\frac{1}{3} k^{2m-1}}{2^m (m-2)!}\left[\frac{\Gamma(m-1) \Gamma(1)}{\Gamma(m-1+1)}(1-2m-2)u(m-2)+2\delta(m-2)\right]\\
\nn &&+\sum_{n=0}^{m-1}\frac{\frac{1}{3}(2n-2m-1)(-1)^n j^{2n}k^{2m-1-2n}}{2^{m}(m-1-n)!n!}\\
\nn &&-\frac{\frac{1}{3}(-2m-1)k^{2m-1}}{2^m(m-1)!}u(m-1)+r_{m,2m-3,2m-2,2m-2}(j,k)\\
\label{Pk3_result} &=&\sum_{n=0}^{m-1}\frac{\frac{1}{3} (2n-2m-1)(-1)^n j^{2n}k^{2m-1-2n}}{2^{m}(m-1-n)!n!}\\
\nn &&+\frac{\frac{2}{3} k^{2m-1}}{2^m (m-2)!}\delta(m-2)+\frac{\frac{1}{3}(2m+1)k^{2m-1}}{2^m(m-1)!}\delta(m-1)\\
\nn &&+r_{m,2m-3,2m-2,2m-2}(j,k).
\ena

Next we find from (\ref{Pkm_parts}) with $t=4$
\bea
\nn \lefteqn{\sum_{q=1}^{k-m} q P_{k,m-1;4}(q+1)-\sum_{q=1}^{j-1} q P_{k,m-1;4}(q+1)} \\
\nn &=&\sum_{q=1}^{k-m} q \left(r_{m,2m-4,2m-4,2m-4}(q+1,k)  \right) \\
\nn &&-\sum_{q=1}^{j-1} q \left(r_{m,2m-4,2m-4,2m-4}(q+1,k) \right)\\
\nn &=&r_{m,2m-2}(k)+r_{m,2m-2,2m-4,2m-2}(j,k)\\
\label{Pk4_result} &=&r_{m,2m-2,2m-2,2m-2}(j,k).
\ena

Combining (\ref{Pk1_result}), (\ref{Pk2_result}), (\ref{Pk3_result}) and (\ref{Pk4_result}) we get
\bea
\nn \lefteqn{\sum_{q=1}^{k-m} q P_{k,m-1}(q+1)-\sum_{q=1}^{j-1} q P_{k,m-1}(q+1)} \\
\label{big_result1} &=&\sum_{n=0}^{m} \frac{(-1)^n j^{2n} k^{2m-2n}}{2^m (m-n)!n!}u(m-1)\\
\nn &&-\frac{2^{m-1}(m-1)! k^{2m-1}}{(2m-1)!}u(m-1)+\frac{k^{2m-1}}{2^m (m-1)!}(3-2m)\delta(m-1)\\ 
\nn &&-\sum_{n=0}^{m-1} \frac{(-1)^n j^{2n+1} k^{2m-2-2n}}{2^m (m-1-n)!n!} \frac{(2n-1)}{(2n+1)}+r_{m,2m-2,2m-2,2m-2}(j,k)\\
\nn &&+\frac{2^{m-1}(m-1)! k^{2m-1}}{(2m-1)!}u(m-2)-\frac{\frac{1}{3} k^{2m-1}}{2^{m-1}(m-2)!}\delta(m-2)\\
\nn &&-\sum_{n=0}^{m-2} \frac{\frac{1}{3}(2n-3)(-1)^{n+1} j^{2n+3} k^{2m-4-2n}}{2^{m-1}(m-2-n)!n!(2n+3)}\\
\nn &&+r_{m,2m-2,2m-2,2m-2}(j,k)\\
\nn &&+\sum_{n=0}^{m-1}\frac{\frac{1}{3}(2n-2m-1)(-1)^n j^{2n}k^{2m-1-2n}}{2^{m}(m-1-n)!n!}\\
\nn &&+\frac{\frac{2}{3} k^{2m-1}}{2^m (m-2)!}\delta(m-2)+\frac{\frac{1}{3}(2m+1)k^{2m-1}}{2^m(m-1)!}\delta(m-1)\\
\nn &&+r_{m,2m-3,2m-2,2m-2}(j,k)+r_{m,2m-2,2m-2,2m-2}(j,k).
\ena

Note that for all integers $m$
\bea
\nn 0&=& -\frac{2^{m-1}(m-1)! k^{2m-1}}{(2m-1)!}u(m-1)+\frac{k^{2m-1}}{2^m (m-1)!}(3-2m)\delta(m-1)\\ 
\nn &&+\frac{2^{m-1}(m-1)!k^{2m-1}}{(2m-1)!}u(m-2)-\frac{\frac{1}{3} k^{2m-1}}{2^{m-1}(m-2)!}\delta(m-2)\\
\nn &&+\frac{\frac{2}{3} k^{2m-1}}{2^m (m-2)!}\delta(m-2)+\frac{\frac{1}{3}(2m+1)k^{2m-1}}{2^m(m-1)!}\delta(m-1)
\ena
and
\beas
\nn &-&\sum_{n=0}^{m-1} \frac{(-1)^n j^{2n+1} k^{2m-2-2n}}{2^m (m-1-n)!n!} \frac{(2n-1)}{(2n+1)}-\sum_{n=0}^{m-2} \frac{\frac{1}{3}(2n-3)(-1)^{n+1} j^{2n+3} k^{2m-4-2n}}{2^{m-1}(m-2-n)!n!(2n+3)}\\
\nn &=&\sum_{n=0}^{m-1} \frac{(-1)^{n+1} \frac{1}{3}(2n-3)j^{2n+1} k^{2m-2-2n}}{2^m (m-1-n)!n!}.\\
\enas
so noting that $m-1\geq 0$ in this induction proof (\ref{big_result1}) becomes
\bea
\nn \lefteqn{\sum_{q=1}^{k-m} q P_{k,m-1}(q+1)-\sum_{q=1}^{j-1} q P_{k,m-1}(q+1)} \\
\nn &=&\sum_{n=0}^{m} \frac{(-1)^n j^{2n} k^{2m-2n}}{2^m (m-n)!n!}\\
\nn &&+\sum_{n=0}^{m-1} \frac{\frac{1}{3}(2n-3)(-1)^{n+1} j^{2n+1} k^{2m-2-2n}}{2^{m-1}(m-1-n)!n!}\\
\nn &&+\sum_{n=0}^{m-1}\frac{\frac{1}{3}(2n-2m-1)(-1)^n j^{2n}k^{2m-1-2n}}{2^{m}(m-1-n)!n!}\\
\nn &&+r_{m,2m-2,2m-2,2m-2}(j,k)
\ena
which proves (\ref{P_is_P_tilde_plus_r}).

Now we may write (\ref{P_is_P_tilde_plus_r}) as
\beas
\nn \lefteqn{P_{k,m}(j)} \\
&=&\frac{(-1)^m j^{2m}}{2^m m!} + \sum_{n=0}^{m-1} \frac{(-1)^{n} j^{2n} k^{2m-2n}}{2^m (m-n)! n!}\\
\nn &&+\frac{\frac{1}{3}(2m-5)(-1)^{m} j^{2m-1}}{2^m(m-1)!}+\sum_{n=0}^{m-2} \frac{\frac{1}{3}(2n-3)(-1)^{n+1} j^{2n+1} k^{2m-2-2n}}{2^m(m-1-n)!n!}\\
\nn &&+\sum_{n=0}^{m-1} \frac{\frac{1}{3}(2n-2m-1)(-1)^n  j^{2n} k^{2m-1-2n}}{2^m(m-1-n)!n!}+r_{m,2m-2,2m-2,2m-2}(j,k)\\
&=&(-1)^m\left(\frac{j^{2m}+\frac{1}{3}m(2m-5)j^{2m-1}}{2^m m!}\right)+r_{(k,m),2m-2}(j)
\enas
which proves (\ref{P-km-1}). Now interchanging $k$ and $j$ and letting $m=v$ (\ref{P_is_P_tilde_plus_r}) becomes
\beas
\nn \lefteqn{P_{j,v}(k)} \\
&=&\sum_{n=0}^{v} \frac{(-1)^{n}  k^{2n} j^{2v-2n}}{2^v (v-n)! n!}+\sum_{n=0}^{v-1} \frac{\frac{1}{3}(2n-3)(-1)^{n+1} k^{2n+1} j^{2v-2-2n}}{2^v(v-1-n)!n!} \\
\nn &&+\sum_{n=0}^{v-1} \frac{\frac{1}{3}(2n-2v-1)(-1)^n  k^{2n} j^{2v-1-2n}}{2^v(v-1-n)!n!}+r_{v,2v-2,2v-2,2v-2}(k,j)
\enas
so
\beas
P_{j,v}(1)=\frac{j^{2v}-\frac{1}{3}v(2v+1)j^{2v-1}}{2^v v!}+r_{v,2v-2}(j)
\enas
which proves (\ref{P_{j,v}}). The proof of the lemma is now complete.
\bbox

\begin{lemma} With $P_{j,v}^0$ given in Definition \ref{P_def},
\bea \label{fN_j}
(f_NN)_j=\sum_{v=0}^j (-1)^v P_{j,v}^0(1)(f_NN)^{j-v} \quad \mbox{for $j \ge 0$,}
\ena
and
\bea \label{N-j_k-j}
(N-k)_{j-k}=\sum_{v=0}^{j-k} (-1)^v P_{j,v}^0(k)N^{j-k-v} \quad \mbox{for $j -k \ge 0$.}
\ena
\end{lemma}
\noindent {\bf Proof:} First note that the validity of (\ref{N-j_k-j}) establishes (\ref{fN_j}) since with $k=0$, $P_{j,v}^0(0)=P_{j,v}^0(1)$.   Expanding $(x-k)_j$ in $x \in \mathbb{R}$
for $k \in \mathbb{N}$ and $j \ge 0$ we obtain
\beas
(x-k)_j=\sum_{v=0}^j (-1)^v a_{k+j-1,v}(k) x^{j-v},
\enas
where
\beas
a_{j,v}(u)=\sum_{u \leq l_1<l_2<\cdots <l_v\leq j} \prod_{m=1}^v l_m
\enas
with any empty product set to 1. Comparison with (\ref{fN_j}) shows that it suffices to prove
\bea
\label{av=P}
a_{k+j-1,v}(k)=P_{j,v}^0(k).
\ena
Equality holds for $v=0$. Assuming (\ref{av=P}) holds for some $v-1 \ge 0$, we have
\beas a_{k+j-1,v}(k)
&=& \sum_{k\leq l_1<l_2<\cdots <l_v\leq k+j-1} \prod_{m=1}^v l_m \\
&=& \sum_{k\leq q \leq k+j-v} q \sum_{l_2,\ldots,l_v: q<l_2<\cdots <l_v \leq k+j-1}
\prod_{m=2}^v l_m \\
&=& \sum_{k \leq q \leq k+j-v} q \sum_{l_1,\ldots,l_{v-1}:q+1\leq l_1<l_2<\cdots <l_{v-1} \leq k+j-1}
\prod_{m=1}^{v-1} l_m \\
&=& \sum_{q=k}^{k+j-v} q a_{k+j-1,v-1}(q+1).
\enas
As  $P_{j,v}^0(k)$ is characterized by the same recursion, the inductive step is complete, proving (\ref{N-j_k-j}).  \bbox
\clearpage
\noindent {\bf Proof of Theorem 1.1}

Expanding out the product in (\ref{def-Corr-H}), by
(\ref{E-prod-1A}) and (\ref{def-fN}) we have \bea \nonumber
\mbox{Corr}(k)&=& E\left( \prod_{A \in H} (I_A-f_N)
\right)\\
&=& \nonumber E  \sum_{G \subset H} \left( \prod_{A \in G}I_A
\right)
(-f_N)^{|H \setminus G|}\\
&=& \nonumber \sum_{G \subset H} \frac{{N-|G| \choose n-|G|}}{{N
\choose n}}
(-f_N)^{|H|-|G|}\\
&=& \nonumber \sum_{j=0}^k \sum_{G \subset H, |G|=j} \frac{{N-j
\choose n-j}}{{N \choose n}}
(-f_N)^{k-j}\\
&=& \nonumber \sum_{j=0}^{k} \frac{{k\choose j}{{N-j\choose
n-j}}}{{N\choose
n}}\left(-f_N\right)^{k-j}\\
&=& \nonumber \sum_{j=0}^k {k\choose
j}\frac{(n)_j}{(N)_j}(-f_N)^{k-j}\\
&=& \nonumber \frac{\sum_{j=0}^{k} {k\choose j} (n)_j (N-j)_{k-j}
(-f_N)^{k-j}}{(N)_k} \\
\nonumber &=&\frac{\alpha(k,f_N)}{(N)_k}, \ena where \bea
\label{def-alpha} \alpha(k,f_N) =(-f_N)^k \sum_{j=0}^{k} {k\choose
j}(-1)^j \lambda(k,j,f_N), \quad \mbox{for $f_N \in (0,1)$}\ena and
\bea
\label{def-lambdak}
\lambda(k,j,f_N)=f_N^{-j}(f_NN)_j (N-j)_{k-j}, \quad \mbox{for
$j=0,1,2,\ldots,k$}.
\ena
We also write \bea \label{corrkf}
\mbox{Corr}(k,f_N)=\frac{\alpha(k,f_N)}{(N)_k} \quad \mbox{for $f_N \in
(0,1)$.} \ena Using (\ref{fN_j}) and (\ref{N-j_k-j}), and recalling that $P_{j,v}^0(1)=0$ when $v <0$ then (\ref{def-lambdak}) becomes 
\beas
\lambda(k,j,f_N)&=& f_N^{-j}\left(\sum_{v=0}^j (-1)^v
P_{j,v}^0(1)(f_NN)^{j-v}\right)
\left(\sum_{i=0}^{k-j} (-1)^i P_{k,i}^0(j) N^{k-j-i}\right)\\
&=& f_N^{-j} \sum_{r=0}^k \left(  \sum_{v=j-r}^{k-r} (-1)^v P_{j,v}^0(1)f_N^{j-v} (-1)^{k-v-r} P_{k,k-v-r}^0(j) \right) N^r\\
&=& \sum_{r=0}^k  \sum_{v=j-r}^{k-r}  f_N^{-v} P_{j,v}^0(1) (-1)^{k-r} P_{k,k-v-r}^0(j) N^r\\
&=& \sum_{v=j-k}^k f_N^{-v} P_{j,v}^0(1) \sum_{r=j-v}^{k-v}  (-1)^{k-r} P_{k,k-v-r}^0(j) N^r\\
&=& \sum_{v=0}^k f_N^{-v} P_{j,v}^0(1) \sum_{r=j-v}^{k-v}  (-1)^{k-r} P_{k,k-v-r}^0(j) N^r,
\enas
since $P_{j,v}^0(1)=0$ for $v<0$.

For a function represented as the power series $\mu(f)=\sum_{j=-\infty}^\infty c_j f^j$, for $k \in \mathbb{Z}$ let $\mu(f;k)$ denote the coefficient of $f^k$ in
$\mu(f)$, that is,
$$
\mu(f;k)=c_k.
$$

Making the change of variable $m=k-r$ in the first sum below, and again
using that $P_{k,m}^0(j)=0$ for $m<0$, for $0 \le v \le k$ we obtain\\

\bea \nonumber \lambda(k,j,f_N;-v) &=&
P_{j,v}^0(1)\sum_{r=j-v}^{k-v} (-1)^{k-r} P_{k,k-v-r}^0(j)
N^r\\
\nonumber &=& P_{j,v}^0(1)\sum_{m=k-j+v}^v (-1)^m P_{k,m-v}^0(j)N^{k-m} \\
\nonumber &=& P_{j,v}^0(1)\sum_{m=v}^{k-j+v} (-1)^m P_{k,m-v}^0(j)N^{k-m} \\
\label{using_P_a} &=& \sum_{m=v}^{k-j+v}\lambda(k,j,f_N;-v:m) \ena where
\beas
\lambda(k,j,f_N;-v:m)= (-1)^m P_{j,v}^0(1) P_{k,m-v}^0(j)N^{k-m}.
\enas

Now, by (\ref{def-alpha}) and (\ref{using_P_a}), \bea \nn
\alpha(k,f_N;k-v) &=&(-1)^k \sum_{j=0}^{k} {k\choose j}(-1)^j
\lambda(k,j,f_N;-v)\\
\nn &=&(-1)^k \sum_{j=0}^{k}
 \sum_{m=v}^{k-j+v} {k\choose j}(-1)^j \lambda(k,j,f_N;-v:m)\\
\nn &=&(-1)^k  \sum_{m=v}^{k+v}
  \sum_{j=0}^{k-m+v} {k\choose j}(-1)^j \lambda(k,j,f_N;-v:m)\\
\nn &=&(-1)^k  \sum_{m=v}^{k+v}
  \sum_{j=0}^{k-m+v} {k\choose j}(-1)^{j+m} P_{j,v}^0(1)
P_{k,m-v}^0(j)N^{k-m}\\
\nn &=&(-1)^k  \sum_{m=v}^{k+v}
  \sum_{j=0}^k {k\choose j}(-1)^{j+m} P_{j,v}^0(1)
P_{k,m-v}^0(j)N^{k-m}\\
\label{alpha_result} &=&(-1)^k  \sum_{m=v}^{k+v}
  \sum_{j=0}^k {k\choose j}(-1)^{j+m} P_{j,v}(1)
P_{k,m-v}(j)N^{k-m},
\ena
where in the last two equalities we invoke Lemma \ref{P-lemma} to apply
$P_{k,m-v}^0(j)=0$ for $j>k-m+v$, and
$P_{j,v}^0(1)=P_{j,v}(1)$, $P_{k,m-v}^0(j)=P_{k,m-v}(j)$ for $0\leq j\leq k$.

At this point in our proof we will consider $k$ even and $k$ odd separately. For $k$ even it will be convenient to write (\ref{P-km-1}) and (\ref{P_{j,v}}), respectively, as
\bea 
\nn P_{k,m}(j)=(-1)^m\frac{j^{2m}}{2^m m!}+r_{(k,m),2m-1}(j),
\ena
and 
\bea 
\nn P_{j,v}(1)=\frac{j^{2v}}{2^v v!}+r_{v,2v-1}(j).
\ena
Substituting these last two expressions into (\ref{alpha_result}) we get
\bea
 \nn
\nn \alpha(k,f_N;k-v)&=& (-1)^k \sum_{m=v}^{k+v}\sum_{j=0}^{k} \Biggl{[}{k\choose j}(-1)^{j+m} \left(\frac{j^{2v}}{2^v v!}+r_{v,2v-1}(j)\right)\\
\nn &&\times\Biggl{(}\frac{(-1)^{m-v} j^{2(m-v)}}{2^{m-v} (m-v)!}+r_{(k,m,v),2(m-v)-1}(j)\Biggl{)}\Biggr{]}N^{k-m}\\
\nn &=& (-1)^k \sum_{m=v}^{k+v} \sum_{j=0}^{k} \Biggl{[}{k\choose j}(-1)^{j+m}\\
\label{alpha_beta} &&\times \Biggl{(}\frac{(-1)^{m-v} j^{2m}}{2^{m} v!(m-v)!}+r_{(k,m,v),2m-1}(j)\Biggl{)}\Biggr{]}N^{k-m}.
\ena

From (\ref{corrkf}) and (\ref{alpha_beta}) we obtain
\bea
\nn \lefteqn{N^{k/2}\mbox{Corr}(k,f_N;k-v)} \\
\nn &=&N^{k/2}\frac{\alpha(k,f_N;k-v)}{(N)_k} \\
\nn &=&\frac{N^{k/2}}{(N)_k}(-1)^k \sum_{m=v}^{k+v} \sum_{j=0}^{k} \Biggl{[}{k\choose j}(-1)^{j+m}\\
\nn && \times \Biggl{(}\frac{(-1)^{m-v} j^{2m}}{2^{m} v!(m-v)!}+r_{(k,m,v),2m-1}(j)\Biggl{)}\Biggr{]}N^{k-m}
\ena
which becomes
\bea
\nn \lefteqn{N^{k/2}\mbox{Corr}(k,f_N;k-v)} \\
\label{Corr_result}  &=&\frac{N^{k/2}}{(N)_k}\frac{(-1)^{k-v}}{v!} \sum_{m=v}^{k+v} \left[\frac{1}{2^{m} (m-v)!}
\sum_{j=0}^{k} {k\choose j}(-1)^{j} j^{2m}\right.  \\
\nn &&\left.+\sum_{j=0}^{k} {k\choose j}(-1)^{j} r_{(k,m,v),2m-1}(j)\right]N^{k-m}.
\ena

Using (\ref{Stirling 2}), (\ref{Corr_result}) becomes
\bea \label{nonlim_corr_result}
\nn \lefteqn{N^{k/2}\mbox{Corr}(k,f_N;k-v)} \\
&=&\frac{N^{k/2}}{(N)_k}\frac{(-1)^{k-v}}{v!} \Biggl{(}\sum_{m=v}^{k+v} \frac{1}{2^{m} (m-v)!} (-1)^k k!\left\{2m \atop k \right\}  \\
\nn &&+\sum_{m=v}^{k+v} \sum_{j=0}^{k} {k\choose j}(-1)^{j} r_{(k,m,v),2m-1}(j)\Biggl{)}N^{k-m}. 
\ena

\noindent Since $(N)_k$ is of degree $k$, then when $m>k/2$,
\beas
\lim_{N \rightarrow \infty} \frac{N^{k/2}}{(N)_k} N^{k-m} =\lim_{N \rightarrow \infty} \frac{N^{k}}{(N)_k} \frac{N^{k/2}}{N^m}=0
\enas
\noindent and when $m \le k/2$, from (\ref{Stirling 2}) and (\ref{second_Stirlings}) we have that
\beas
\sum_{j=0}^k {k\choose j}(-1)^j j^p = (-1)^k k! \left\{p \atop k \right\}=0 \quad \mbox{for all $p \le 2m-1$,}
\enas
and therefore, by linearity,
\beas
\sum_{j=0}^{k} {k\choose j}(-1)^{j} r_{(k,m,v),2m-1}(j)=0 \quad \mbox{for $m \le k/2$,}
\enas
and since for $k$ even
\beas
\left\{2m \atop k \right\}=0 \quad \mbox{for $m<k/2$,}
\enas
then upon letting $N\rightarrow \infty$, (\ref{nonlim_corr_result}) is zero for $v>k/2$ and for $v\leq k/2$ \mbox{(\ref{nonlim_corr_result}) becomes}
\bea \label{Corr_even_result}
\lefteqn{\lim_{N \rightarrow \infty} N^{k/2}\mbox{Corr}(k,f_N;k-v)} \\
\nn &=&\lim_{N \rightarrow \infty} \frac{N^{k/2}}{(N)_k}\frac{(-1)^{k-v}}{v!}\left(\sum_{m=k/2}^{k/2} \frac{1}{2^{m} (m-v)!} (-1)^k k!\left\{2m \atop k \right\}\right)N^{k-m}. 
\ena

\noindent When $m=k/2$,
\beas
\lim_{N \rightarrow \infty} \frac{N^{k/2}}{(N)_k} N^{k-m} =\lim_{N \rightarrow \infty} \frac{N^{k}}{(N)_k} \frac{N^{k/2}}{N^m}=\lim_{N \rightarrow \infty} \frac{N^{k}}{(N)_k} \frac{N^{k/2}}{N^{k/2}}=1.
\enas
Therefore, using (\ref{second_Stirlings}), (\ref{Corr_even_result}) becomes
\bea 
\nn \lim_{N \rightarrow \infty} N^{k/2}\cdot \mbox{Corr}(k,f_N;k-v)&=& (-1)^{v}  \frac{k!}{v! 2^{k/2} \left(\frac{k}{2}-v\right)!}\\
\nn &=& (-1)^{v}  \frac{k!}{2^{k/2} \left(\frac{k}{2}\right)!}{k/2 \choose v}\\
\nn &=& (-1)^v  EZ^k {k/2 \choose v}.
\ena
Now, by our previous notation, $\mbox{Corr}(k)=\mbox{Corr}(k,f_N)$, so
\beas
\lim_{N \rightarrow \infty} N^{k/2}\cdot \mbox{Corr}(k)&=&\lim_{N \rightarrow \infty} N^{k/2}\cdot \mbox{Corr}(k,f_N)\\
&=&\lim_{N \rightarrow \infty} N^{k/2}\sum_{v=0}^k \mbox{Corr}(k,f_N;k-v)f_N^{k-v}\\
&=&\sum_{v=0}^{k/2} {k/2 \choose v} f^{k-v} (-1)^v EZ^k\\
&=&[f(f-1)]^{k/2}EZ^k
\enas
which proves the theorem for $k$ even.

For $k$ odd we have from (\ref{P-km-1}), (\ref{P_{j,v}}) and (\ref{alpha_result}) that
\bea
 \nn
\lefteqn{\alpha(k,f_N;k-v)} \\
\nn &=& (-1)^k \sum_{m=v}^{k+v}\sum_{j=0}^{k} {k\choose j}(-1)^{j+m} \left(\frac{j^{2v}-\frac{1}{3}v(2v+1)j^{2v-1}}{2^v v!}+r_{v,2v-2}(j)\right)\\
\nn &&\times\Biggl{[}(-1)^{m-v}\left(\frac{j^{2(m-v)}+\frac{1}{3}(m-v)(2(m-v)-5)j^{2(m-v)-1}}{2^{m-v} (m-v)!}\right)\\
\nn &&+r_{(k,m,v),2(m-v)-2}(j)\Biggl{]}N^{k-m}
\ena
which becomes
\bea
 \nn
\lefteqn{\alpha(k,f_N;k-v)} \\
\nn &=& (-1)^k \sum_{m=v}^{k+v}\sum_{j=0}^{k} {k\choose j}(-1)^{j+m}\\
\nn &&\times\Biggl{[}(-1)^{m-v}\left(\frac{j^{2m}+\frac{1}{3}\left[(m-v)(2m-2v-5)-v(2v+1)\right]j^{2m-1}}{2^m v! (m-v)!}\right)\\
\nn &&+r_{(k,m,v),2m-2}(j)\Biggl{]}N^{k-m}.
\ena

From (\ref{corrkf}) and (\ref{alpha_beta}) we obtain
\bea
\nn \lefteqn{N^{(k+1)/2}\mbox{Corr}(k,f_N;k-v)=N^{(k+1)/2}\frac{\alpha(k,f_N;k-v)}{(N)_k}}\\
\nn &=&\frac{N^{(k+1)/2}}{(N)_k}(-1)^k \sum_{m=v}^{k+v}\sum_{j=0}^{k} {k\choose j}(-1)^{j+m}\\
\nn &&\times\Biggl{[}(-1)^{m-v}\left(\frac{j^{2m}+\frac{1}{3}\left[(m-v)(2m-2v-5)-v(2v+1)\right]j^{2m-1}}{2^m v! (m-v)!}\right)\\
\nn &&+r_{(k,m,v),2m-2}(j)\Biggl{]}N^{k-m}\\
\label{Corr_result2}  &=&\frac{N^{(k+1)/2}}{(N)_k}\frac{(-1)^{k-v}}{v!} \sum_{m=v}^{k+v} \left[\frac{1}{2^{m} (m-v)!}
\sum_{j=0}^{k} {k\choose j}(-1)^{j} \right.\\
\nn && \times\left( j^{2m}+\frac{1}{3}\left[(m-v)(2m-2v-5)-v(2v+1)\right]j^{2m-1}\right)  \\
\nn && \left.+\sum_{j=0}^{k} {k\choose j}(-1)^{j} r_{(k,m,v),2m-2}(j)\right]N^{k-m}.
\ena

Using (\ref{Stirling 2}), (\ref{Corr_result2}) becomes
\bea \label{Corr_result3}
\nn \lefteqn{N^{(k+1)/2}\mbox{Corr}(k,f_N;k-v)}\\
    &=&\frac{N^{(k+1)/2}}{(N)_k}\frac{(-1)^{k-v}}{v!} \left[\sum_{m=v}^{k+v} \frac{1}{2^{m} (m-v)!}(-1)^k k! \left\{2m \atop k \right\}\right.\\
\nn &&+\sum_{m=v}^{k+v} \frac{\frac{1}{3}\left[(m-v)(2m-2v-5)-v(2v+1)\right]}{2^{m} (m-v)!}(-1)^k k! \left\{2m-1 \atop k \right\}  \\
\nn &&+\sum_{m=v}^{k+v} \left. \sum_{j=0}^{k} {k\choose j}(-1)^{j} r_{(k,m,v),2m-2}(j)\right]N^{k-m}.
\ena

\noindent Since $(N)_k$ is of degree $k$, then when $m>(k+1)/2$,
\beas
\lim_{N \rightarrow \infty} \frac{N^{(k+1)/2}}{(N)_k} N^{k-m} =\lim_{N \rightarrow \infty} \frac{N^{k}}{(N)_k} \frac{N^{(k+1)/2}}{N^m}=0
\enas
\noindent and when $m \le (k+1)/2$, from (\ref{Stirling 2}) and (\ref{second_Stirlings}) we have that
\beas
\sum_{j=0}^k {k\choose j}(-1)^j j^p = (-1)^k k! \left\{p \atop k \right\}=0 \quad \mbox{for all $p \le 2m-2$,}
\enas
and therefore, by linearity,
\beas
\sum_{j=0}^{k} {k\choose j}(-1)^{j} r_{(k,m,v),2m-2}(j)=0 \quad \mbox{for $m \le (k+1)/2$},
\enas

\noindent and since for $k$ odd
\beas
\left\{2m \atop k \right\}=0 \quad \mbox{for $m<(k+1)/2$},
\enas
and
\beas
\left\{2m-1 \atop k \right\}=0 \quad \mbox{for $m<(k+1)/2$},
\enas
then upon letting $N\rightarrow \infty$ (\ref{Corr_result3}) is zero for $v>(k+1)/2$ and for $v\leq (k+1)/2$ (\ref{Corr_result3}) becomes
%
\bea
\nn \lefteqn{\lim_{N \rightarrow \infty} N^{(k+1)/2}\mbox{Corr}(k,f_N;k-v)}\\
\label{Corr_odd_result} &=&\lim_{N \rightarrow \infty} \frac{N^{(k+1)/2}}{(N)_k}\frac{(-1)^{k-v}}{v!} \left[\sum_{m=(k+1)/2}^{(k+1)/2} \frac{1}{2^{m} (m-v)!}(-1)^k k! \left\{2m \atop k \right\}\right.\\
\nn &&\left. +\sum_{m=(k+1)/2}^{(k+1)/2}\frac{\frac{1}{3}\left[(m-v)(2m-2v-5)-v(2v+1)\right]}{2^{m} (m-v)!}(-1)^k k! \left\{2m-1 \atop k \right\}\right]N^{k-m}  
\ena

\noindent When $m=(k+1)/2$,
\beas
\lim_{N \rightarrow \infty} \frac{N^{(k+1)/2}}{(N)_k} N^{k-m} =\lim_{N \rightarrow \infty} \frac{N^{k}}{(N)_k} \frac{N^{(k+1)/2}}{N^m}=\lim_{N \rightarrow \infty} \frac{N^{k}}{(N)_k} \frac{N^{(k+1)/2}}{N^{(k+1)/2}}=1.
\enas
Therefore,(\ref{Corr_odd_result}) becomes
\beas \lefteqn{\lim_{N \rightarrow \infty} N^{(k+1)/2}\cdot \mbox{Corr}(k,f_N;k-v)}\\
\nn &=&\frac{(-1)^{k-v}}{v!} \left[ \frac{1}{2^{(k+1)/2} \left(\frac{k+1}{2}-v\right)!}(-1)^k k! \left\{k+1 \atop k \right\}\right.\\
\nn &&+\left. \frac{\frac{1}{3}\left[\left(\frac{k+1}{2}-v\right)(k+1-2v-5)-v(2v+1)\right]}{2^{(k+1)/2} \left(\frac{k+1}{2}-v\right)!}(-1)^k k! \left\{k \atop k \right\}\right].
\enas
Using (\ref{second_Stirlings}) this last result becomes
\beas
\nn \lefteqn{\lim_{N \rightarrow \infty} N^{(k+1)/2}\cdot \mbox{Corr}(k,f_N;k-v)}\\
\nn &=&\frac{\frac{1}{3}(-1)^{v}k!}{v!2^{(k+1)/2}\left(\frac{k+1}{2}-v\right)!} \left[ 3 {k+1\choose 2}\right.\\
\nn &&+\left. \left(\frac{k+1}{2}-v\right)(k-2v-4)-v(2v+1) \right]\\
\nn &=&\frac{\frac{2}{3}(-1)^{v}(k-1)(k+1-v)k!}{v!2^{(k+1)/2}\left(\frac{k+1}{2}-v\right)!}.
\enas
We may write this last result as
\beas
\nn \lefteqn{\lim_{N \rightarrow \infty} N^{(k+1)/2}\cdot \mbox{Corr}(k,f_N;k-v)}\\
\nn &=&\frac{\frac{2}{3}(-1)^{v}(k-1)(k+1-v)k! {\frac{k+1}{2}\choose v}}{2^{(k+1)/2}\left(\frac{k+1}{2}\right)!}\\
\nn &=&\frac{\frac{2}{3}(-1)^{v}(k-1)(k+1-v) EZ^{k+1}}{k+1} {\frac{k+1}{2}\choose v}.
\enas
Hence,
\beas
\lim_{N \rightarrow \infty} N^{(k+1)/2}\cdot \mbox{Corr}(k)&=&\lim_{N \rightarrow \infty} N^{(k+1)/2}\cdot \mbox{Corr}(k,f_N)\\
&=&\lim_{N \rightarrow \infty} N^{(k+1)/2}\sum_{v=0}^k \mbox{Corr}(k,f_N;k-v)f_N^{k-v}\\
&=&\sum_{v=0}^{k} \frac{\frac{2}{3}(-1)^{v}(k-1)(k+1-v)}{k+1} {\frac{k+1}{2}\choose v} f^{k-v} EZ^{k+1}\\
&=&\frac{\frac{2}{3}(k-1)}{k+1} EZ^{k+1} \frac{d}{df}\left(\sum_{v=0}^{(k+1)/2}  (-1)^{v}{\frac{k+1}{2}\choose v} f^{k+1-v}\right)\\
&=&\frac{\frac{2}{3}(k-1)}{k+1} EZ^{k+1} \frac{d}{df}\left(\left[f(f-1)\right]^{(k+1)/2}\right)\\
&=&[f(f-1)]^{(k-1)/2}(2f-1)\frac{1}{3}(k-1)EZ^{k+1}
\enas
which proves the theorem for $k$ odd.

The proof of the theorem is now complete. \bbox
\\
\section*{Acknowledgments}
I would like to express my sincere gratitude to Professor Larry Goldstein of the Department of Mathematics at the University of Southern California for suggesting this topic and for his many valuable comments that greatly improved this correspondence.
\\
\\
\section*{References}
\begin{enumerate}

\item \label{AS} Abramowitz, M., and Stegun, I., (1965). Handbook of Mathematical Functions.
Ninth edition. Dover, New York, NY.

\item \label{AGL} Arraita, R., Goldstein, L., and Langholz, B.
(2005). Local central limit theorems, the high-order correlations
of rejective sampling and logistic likelihood asymptotics. {\em
Ann. Statist.}, {\bf 33}, pp. 871--914.

\item \label{Garrappa} Garrappa, R. (2007) Some Formulas for Sums of Binomial Coefficients and Gamma Functions. 
{\em International Mathematical Forum}, {\bf 15},  pp. 725-733.

\item \label{Concrete} Graham, R., Knuth, D., and Patashnik, O.
(1994). Concrete Mathematics: A Foundation for Computer Science.
Second edition. Addison-Wesley Publishing Company, Reading, MA.

\item \label{Hajek2} H\'ajek, J. (1964) Asymptotic theory of
rejective sampling with varying probabilities from a finite
population. {\em Ann. Math. Statist.}, {\bf 35},  pp. 1491--1523.

\item \label{Hajek} H\'ajek, J. (1981) Sampling from a finite
population. Edited by Václav Dupa\v c. With a foreword by P. K.
Sen. Statistics: Textbooks and Monographs, 37. Marcel Dekker,
Inc., New York.


\item \label{oeis} The On-Line Encyclopedia of Integer Sequences: http://oeis.org

\end{enumerate}
\end{document}